\documentclass[11pt]{amsart}
\newtheorem{theorem}{Theorem}[section]
\newtheorem{lemma}[theorem]{Lemma}

\theoremstyle{definition}
\newtheorem{definition}[theorem]{Definition}

\theoremstyle{remark}

\numberwithin{equation}{section}
                                                                             
\newcommand{\Tr}{{\text{Tr}}}

\begin{document}


\title[Approximating L$^2$ Betti numbers of amenable covering spaces]
{Approximating L$^2$ invariants of amenable covering spaces: 
A combinatorial approach.}
\author{Jozef Dodziuk}
\address{Department of Mathematics, CUNY, NYC, USA}
\email{jzdqc@cunyvm.cuny.edu}
\author{Varghese Mathai}
\address{Department of Mathematics, University of Adelaide, Adelaide 5005,
Australia}
\email{vmathai@maths.adelaide.edu.au}
\date{SEPTEMBER 1996}
\subjclass{Primary: 58G11, 58G18 and 58G25.}
\keywords{L$^2$ Betti numbers, approximation theorems, amenable groups}

\begin{abstract}
In this paper, we prove that the $L^2$ Betti numbers of an amenable
covering space can be approximated by the average Betti numbers of a 
regular exhaustion, proving a conjecture in \cite{DM}.
We also prove that an arbitrary amenable covering space 
of a finite simplicial complex is of determinant class.  
\end{abstract}
\maketitle

\section*{Introduction}

Let $Y$ be a connected simplicial complex.  Suppose that
$\pi$ acts freely and simplicially on $Y$ so that $X=Y/\pi$ is a finite
simplicial complex.  Let ${\mathcal F}$ a finite subcomplex of $Y$, which 
is a fundamental domain
for the action of $\pi$ on $Y$.

We assume that $\pi$ is amenable.  The F\o{}lner criterion for amenability
of $\pi$ enables one to get, cf.  \cite{Ad}, a {\em regular exhaustion} 
$\big\{Y_{m}\big\}^{\infty}_
{m=1}$, that is a sequence of finite subcomplexes of $Y$  
such that

(1) $Y_{m}$ consists of $N_{m}$ translates $g.{\mathcal F}$ of 
${\mathcal F}$ for
$g\in\pi$.

(2) $\displaystyle Y=\bigcup^{\infty}_{m=1}Y_{m}\;$.

(3) If $\dot{N}_{m,\delta}$ denotes the number of translates of 
${\mathcal F}$ 
which have distance (with respect to the word metric in $\pi$) less than or
equal to $\delta$ from a translate of ${\mathcal F}$ having a
non-empty intersection with  the topological boundary $\partial{Y}_{m}$ of 
$Y_{m}$ (we identify here $g.{\mathcal F}$ with
$g$)
then, for every $\delta > 0$,
$$
\lim_{m\rightarrow\infty}\;\frac{{\dot{N}}_{m,\delta}}{N_{m}}=0.
$$

One of our main results is

\begin{theorem}[Amenable Approximation Theorem] 
$\;$Let $Y$ be a connected simplicial complex.
Suppose that $\pi$ is amenable and that $\pi$ acts freely and simplicially
on $Y$ so that $X=Y/\pi$ is a finite simplicial complex.  
Let $\big\{Y_{m}\big\}^{\infty}_
{m=1}$ be a regular exhaustion of $Y$.  Then
$$
\lim_{m\rightarrow\infty}\;\frac{b^{j}(Y_{m})}{N_{m}}=b_{(2)}^{j}(Y:\pi)
\;\;\mbox{ for all }\;\;j\ge 0.
$$
$$
\lim_{m\rightarrow\infty}\;\frac{b^{j}(Y_{m}, 
\partial Y_{m})}{N_{m}}=b_{(2)}^{j}(Y:\pi)
\;\;\mbox{ for all }\;\;j\ge 0.
$$
\end{theorem}

Here $b^{j}(Y_{m})$ denotes the $j^{th}$ Betti number of $Y_m$, 
$b^{j}(Y_{m}, \partial Y_{m})$ denotes the $j^{th}$ Betti number of 
$Y_m$ relative its boundary $\partial Y_m$ and $b_{(2)}^{j}(Y:\pi)$
denotes the $j$th $L^2$ Betti number of $Y$. 
(See the next section for the definition 
of the $L^2$ Betti numbers of a manifold)

\noindent{\bf Remarks.} This theorem proves the main conjecture in the 
introduction of 
an earlier paper \cite{DM}. The combinatorial techniques of this paper 
contrasts
with the heat kernel approach used in \cite{DM}.
Under the assumption dim $H^{k}(Y)<\infty$, a special case of the 
Amenable Approximation Theorem above 
is obtained by combining proofs of
Eckmann \cite{Ec} and Cheeger-Gromov \cite{CG}.
The assumption dim $H^{k}(Y)<\infty$ is very restrictive and essentially 
says
that $Y$ is a fiber bundle over a $B\pi$ with fiber a space with
finite fundamental group. Cheeger-Gromov use this to 
show that the Euler
characteristic of a finite $B\pi$, where $\pi$ contains an infinite
amenable normal subgroup, is zero. Eckmann proved the same result 
in the special
case when $\pi$ itself is an infinite amenable group. 

There is a standing conjecture that any normal
covering space of a finite simplicial complex is of determinant class
(cf. section 4 for the definition of determinant class and for a 
more detailed discussion of what follows).
Let $M$ be a smooth compact manifold, and $X$ triangulation of $M$.
Let $\widetilde M$ be any normal covering space of $M$, and $Y$ be the
triangulation of $\widetilde M$ which projects down to $X$.
Then on $\widetilde M$, there are two notions of determinant class,
one analytic and the other combinatorial. Using results of 
Efremov \cite{E}, Gromov and Shubin
\cite{GS}, one observes as in \cite{BFKM} that the 
combinatorial and analytic 
notions of determinant class coincide.
It was proved in \cite{BFK}) using estimates of L\"uck \cite{L} that any 
{\em residually finite} normal
covering space of a finite simplicial complex is of determinant class, which
gave evidence supporting the conjecture. 
Our interest in this conjecture stems from work on $L^2$ 
torsion \cite{CFM}, \cite{BFKM}. The $L^2$ torsion is a well defined 
element in the determinant line of the reduced $L^2$ cohomology, whenever
the covering space is of determinant class. Our next main theorem
says that any {\em amenable} normal
covering space of a finite simplicial complex is of determinant class,
which gives further evidence supporting the conjecture.

\begin{theorem}[Determinant Class Theorem] 
$\;$Let $Y$ be a connected simplicial complex.
Suppose that $\pi$ is amenable and that $\pi$ acts freely and simplicially
on $Y$ so that $X=Y/\pi$ is a finite simplicial complex. Then $Y$ 
is of determinant class.
\end{theorem}

The paper is organized as follows. In the first section, some 
preliminaries on $L^2$ 
cohomology and amenable groups are
presented. In section 2, the main abstract approximation theorem is 
proved. We 
essentially use the combinatorial analogue of the principle of not feeling
the boundary (cf. \cite{DM}) in Lemma 2.3 and a finite dimensional result 
in \cite{L}, to prove this theorem.
Section 3 contains the proof of the Amenable Approximation Theorem and 
some related 
approximation theorems. 
In section 4, we prove that an arbitrary {\em amenable} normal covering 
space of a finite simplicial complex is of determinant class. 

The second author warmly thanks Shmuel Weinberger for some useful 
discussions. This paper has been inspired by L\"uck's work \cite{L} on 
residually finite groups.

\section{Preliminaries}
  
Let $\pi$ be a finitely generated discrete group and ${\mathcal U}(\pi)$ 
be the von Neumann algebra
generated by the action of $\pi$ on $\ell^{2}(\pi)$ via the left regular
representation. It is the weak (or strong) closure of the complex group 
algebra of 
$\pi$, ${\mathbb C}(\pi)$ acting on $\ell^2(\pi)$ by left translation.
The left regular representation is then a unitary representation
$\rho:\pi\rightarrow{\mathcal U}(\pi)$. Let $\Tr_{{\mathcal  U}(\pi)}$ 
be the faithful normal trace on
${\mathcal U}(\pi)$ defined by the inner product 
$\Tr_{{\mathcal  U}(\pi)}(A) \equiv
(A\delta_e,\delta_e)$ for
$A\in{\mathcal U}(\pi)$ and where $\delta_e\in\ell^{2}(\pi)$ is given by 
$\delta_e(e)=1$ and
$\delta_e(g)=0$ for $g\in\pi$ and $g\neq e$. 

Let $Y$ be a simplicial complex, and $|Y|_j$ denote the set of all 
$p$-simplices
in $Y$. Regarding the orientation of simplices, we use the following 
convention. For
each $p$-simplex $\sigma \in |Y|_j$, we identify $\sigma$ with any other 
$p$-simplex
which is obtained by an {\em even} permutation of the vertices in $\sigma$, 
whereas 
we identify $-\sigma$ with any other $p$-simplex
which is obtained by an {\em odd} permutation of the vertices in $\sigma$.
Suppose that
$\pi$ acts freely and simplicially on $Y$ so that $X=Y/\pi$ is a finite
simplicial complex.  Let ${\mathcal F}$ a finite subcomplex of $Y$, which 
is a fundamental domain
for the action of $\pi$ on $Y$. Consider the Hilbert space of square
summable cochains on $Y$, 
$$
C^j_{(2)}(Y) = \Big\{f\in C^j(Y): \sum_{\sigma\; a\; j-simplex}|f(\sigma)|^2
<\infty \Big\}
$$
Since $\pi$ acts freely on $Y$, we see that there is an isomorphism of
Hilbert $\ell^2(\pi)$ modules,
$$
C^j_{(2)}(Y) \cong C^j(X)\otimes\ell^2(\pi)
$$
Here $\pi$ acts trivially on $C^j(X)$
and via the left regular representation on $\ell^2(\pi)$. Let
$$
d_{j}:C^{j}_{(2)}(Y)\rightarrow C^{j+1}_{(2)}(Y)
$$ 
denote the coboundary operator. It is clearly a bounded linear operator. 
Then the (reduced) $L^2$ cohomology groups of $Y$
are defined to be 
$$
H^j_{(2)}(Y) = \frac{\mbox{ker}(d_j)}{\overline{\mbox{im}(d_{j-1})}}.
$$
Let ${d_j}^*$ denote the Hilbert space adjoint of
$d_{j}$.  One defines the combinatorial Laplacian  
$\Delta_{j} : C^{j}_{(2)}(Y) \rightarrow C^{j}_{(2)}(Y)$ as
$\Delta_j = d_{j-1}(d_{j-1})^{*}+(d_{j})^{*}d_{j}$. 

By the Hodge decomposition theorem in this context, there is an isomorphism
of Hilbert $\ell^2(\pi)$ modules,
$$
H^j_{(2)}(Y)\;\; \cong\;\; \mbox{ker} (\Delta_j).
$$
Let $P_j: C^{j}_{(2)}(Y)\rightarrow \mbox{ker} \Delta_j$ denote the 
orthogonal projection to the kernel of the Laplacian. Then the $L^2$ 
Betti numbers $b_{(2)}^j(Y:\pi)$ are defined as
$$
b_{(2)}^j(Y:\pi) = \Tr_{{\mathcal  U}(\pi)}(P_j).
$$ 

In addition, let $\Delta_j^{(m)}$ 
denote the 
Laplacian on the finite dimensional cochain space $C^j(Y_m)$ or 
$C^j(Y_m,\partial Y_m)$.  We do use the same notation for the two Laplacians
since all proofs work equally well for either case.
Let $D_j(\sigma, \tau) = \left< \Delta_j \delta_\sigma, \delta_\tau\right>$
denote the matrix coefficients of the Laplacian $\Delta_j$ and 
${D_j^{(m)}}(\sigma, \tau) = \left< \Delta_j^{(m)} \delta_\sigma, 
\delta_\tau\right>$
denote the matrix coefficients of the Laplacian $\Delta_j^{(m)}$. Let 
$d(\sigma,\tau)$ denote the 
{\em distance} between $\sigma$ and $\tau$ in the natural simplicial 
metric on $Y$, and
$d_m(\sigma,\tau)$ denote the 
{\em distance} between $\sigma$ and $\tau$ in the natural simplicial 
metric on $Y_m$.  This distance (cf. \cite{Elek}) is defined as follows.
Simplexes $\sigma$ and $\tau$ are one step apart, $d(\sigma,\tau)=1$, if 
they have equal dimensions,
$\dim \sigma = \dim \tau =j$, and there exists either a simplex of 
dimension $j-1$ contained in both $\sigma$ and $\tau$ or a simplex of
dimension $j+1$ containing both $\sigma$ and $\tau$.  The distance between
$\sigma$ and $\tau$ is equal to $k$
if there exists a finite sequence $\sigma = \sigma_0, \sigma_1, \ldots , 
\sigma_k = \tau$, $d(\sigma_i,\sigma_{i+1})=1$ for $i=0,\ldots,k-1$, and 
$k$ is the minimal length of such a sequence. 

Then one has the following, which is an easy
generalization of Lemma 2.5 in \cite{Elek} and follows immediately from
the definition of combinatorial Laplacians and finiteness of the 
complex $X=Y/\pi$.  

\begin{lemma} $D_j(\sigma, \tau) = 0$ whenever $d(\sigma,\tau)>1$ and 
$ {D_j^{(m)}}(\sigma, \tau) = 0$ whenever $d_m(\sigma,\tau)>1$. 
There is also a positive constant $C$ independent of $\sigma,\tau$ such that 
$D_j(\sigma, \tau)\le C$ and ${D_j^{(m)}}(\sigma, \tau)\le C$.
\end{lemma}

Let $D_j^k(\sigma, \tau) = \left< \Delta_j^k \delta_\sigma, 
\delta_\tau\right>$
denote the matrix coefficient of the $k$-th power of the Laplacian,  
$\Delta_j^k$, and 
$D_j^{(m)k}(\sigma, \tau) = \left< \left(\Delta_j^{(m)}\right)^k 
\delta_\sigma, \delta_\tau\right>$
denote the matrix coefficient of the $k$-th power of the Laplacian,  
$\Delta_j^{(m)k}$. Then 
$$D_j^k(\sigma, \tau) = \sum_{\sigma_1,\ldots\sigma_{k-1} \in |Y|_j}
D_j(\sigma, \sigma_1)\ldots D_j(\sigma_{k-1}, \tau)$$
and 
$$D_j^{(m)k}(\sigma, \tau) = \sum_{\sigma_1,\ldots\sigma_{k-1} \in |Y_m|_j}
{D_j^{(m)}}(\sigma, \sigma_1)\ldots {D_j^{(m)}}(\sigma_{k-1}, \tau).$$
Then the following lemma follows easily from Lemma 1.1.

\begin{lemma} Let $k$ be a positive integer. Then $D_j^k(\sigma, \tau) = 0$  
whenever $d(\sigma,\tau)>k$ and
$D_j^{(m)k}(\sigma, \tau) = 0$ whenever $d_m(\sigma,\tau)>k$.
There is also a positive constant $C$ independent of $\sigma,\tau$ such that 
$D_j^k(\sigma, \tau)\le C^k$ and $D_j^{(m)k}(\sigma, \tau)\le C^k$.
\end{lemma}

Since $\pi$ commutes with the Laplacian $\Delta_j^k$, it follows that 
\begin{equation}\label{inv}
D_j^k(\gamma\sigma, \gamma\tau) = D_j^k(\sigma, \tau) 
\end{equation}
for all $\gamma\in \pi$ and for all $\sigma, \tau \in |Y|_j$. The 
{\em von Neumann 
trace} of $\Delta_j^k$ is by definition
\begin{equation} \label{vNt}
\Tr_{{\mathcal  U}(\pi)}(\Delta_j^k) = 
\sum_{\sigma\in |X|_j} D_j^k(\sigma, \sigma),
\end{equation}
where $\tilde{\sigma}$ denotes an arbitrarily chosen lift of $\sigma$ to
$Y$.  The trace is well-defined in view of (\ref{inv}).

\subsection{Amenable groups} 
Let $d_1$ be the word metric on $\pi$.  Recall the following
characterization of amenability due
to F\o{}lner, see also \cite{Ad}.

\begin{definition} A discrete group $\pi$ is said to be {\em amenable} 
if there is a sequence
of finite subsets $\big\{\Lambda_{k}\big\}^{\infty}_{k=1}$ such that for
any fixed $\delta>0$
$$
\lim_{k\rightarrow\infty}\;\frac{\#\{\partial_{\delta}\Lambda_{k}\}}{\#
\{\Lambda_{k}\}}=0
$$
where $\partial_{\delta}\Lambda_{k}=
\{\gamma\in\pi:d_1(\gamma,\Lambda_{k})<\delta$
and $d_1(\gamma,\pi-\Lambda_{k})<\delta\}$ is a $\delta$-neighborhood of
the boundary of $\Lambda_{k}$.  Such a sequence 
$\big\{\Lambda_{k}\big\}^{\infty}_
{k=1}$ is called a {\em regular sequence} in $\pi$.  If in addition
$\Lambda_{k}\subset\Lambda_{k+1}$ for all $k\geq 1$ and 
$\displaystyle\bigcup^
{\infty}_{k=1}\Lambda_{k}=\pi$, then the sequence 
$\big\{\Lambda_{k}\big\}^{\infty}_
{k=1}$ is called a {\em regular exhaustion} in $\pi$.
\end{definition}

Examples of amenable groups are:

\begin{itemize}
\item[(1)]Finite groups;
\item[(2)] Abelian groups; 
\item[(3)] nilpotent groups and solvable groups;
\item[(4)] groups of subexponential growth;
\item[(5)] subgroups, quotient groups and extensions of amenable groups;
\item[(6)] the union of an increasing family of amenable groups.

\end{itemize}

Free groups and fundamental groups of closed negatively curved manifolds are
{\em not} amenable.  

Let $\pi$ be a finitely generated amenable discrete group, and
$\big\{\Lambda_{m}\big\}^{\infty}_{m=1}$ a regular exhaustion in $\pi$.
Then it defines a regular exhaustion $\big\{Y_m\big\}^{\infty}_{m=1}$ of
$Y$. 

Let $\{P_j(\lambda):\lambda\in[0,\infty)\}$ denote the right continuous
family of spectral projections of the Laplacian $\Delta_j$.
Since $\Delta_j$ is $\pi$-equivariant, so are $P_j(\lambda) = 
\chi_{[0,\lambda]}(\Delta_j)$, 
for $\lambda\in [0,\infty)$.  Let 
$F:[0,\infty)\rightarrow[0,\infty)$ denote the spectral
density function,
$$
F(\lambda)=\Tr_{{\mathcal U}(\pi)}(P_j(\lambda)).
$$
Observe that the $j$-th $L^{2}$ Betti number of $Y$ is also given by
$$
b_{(2)}^j(Y:\pi)=F(0).
$$
We have the spectral density function for every dimension $j$ but we do not
indicate explicitly this dependence.  All our arguments are performed with a
fixed value of $j$.  

Let $E_{m}(\lambda)$ denote the number of eigenvalues $\mu$ of
$\Delta_j^{(m)}$ satisfying $\mu\leq\lambda$ and which are counted with
multiplicity. We may sometimes omit the subscript $j$ on $\Delta_j^{(m)}$ 
and $\Delta_j$ to simplify the notation.

We next make the following definitions,
$$
\begin{array}{lcl}
F_{m}(\lambda)& = &\displaystyle\frac{E_{m}(\lambda)}{N_{m}}\;\\[+10pt]
\overline{F}(\lambda) & = & \displaystyle\limsup_{m\rightarrow\infty}
F_{m}(\lambda) \;\\[+10pt]
\mbox{\underline{$F$}}(\lambda) & = & \displaystyle\liminf_{m\rightarrow
\infty}F_{m}(\lambda)\; \\[+10pt]
\overline{F}^{+}(\lambda) & = & \displaystyle\lim_{\delta\rightarrow +0}
\overline{F}(\lambda+\delta) \;\\[+10pt]
\mbox{\underline{$F$}}^{+}(\lambda) & = & \displaystyle\lim_{\delta
\rightarrow +0}\mbox{\underline{$F$}}(\lambda+\delta).
\end{array}
$$

\section{Main Technical Theorem}

Our main technical result is

\begin{theorem} Let $\pi$ be countable, amenable group. 
In the notation of section 1, one has
\begin{itemize}
\item[(1)]$\;F(\lambda)=\overline{F}^{+}(\lambda)=
\mbox{\underline{$F$}}^{+}(\lambda)$.
\item[(2)] $\;\overline{F}$ and {\underline{$F$}} are right
continuous at zero and we have the equalities
\begin{align*}
\overline{F}(0) & =\overline{F}^{+}(0)=F(0)=\mbox{\underline{$F$}}
(0)=\mbox{\underline{$F$}}^{+}(0) \\
\displaystyle & =\lim_{m\rightarrow\infty}F_{m}(0)=\lim_{m\rightarrow
\infty}\;\frac{\#\{E_{m}(0)\}}{N_{m}}\;.
\end{align*}
\item[(3)] $\;$Suppose that $0<\lambda<1$.  Then there is a 
constant $K>1$ such
that
$$
F(\lambda)-F(0)\leq-a\;\frac{\log K^{2}}{\log\lambda}\;.
$$
\end{itemize}
\end{theorem}

To prove this Theorem, we will first prove a number of preliminary lemmas.

\begin{lemma} There exists a positive number $K$  such that the operator
norms of 
 $\Delta_j$ and of $\Delta_j^{(m)}$ for all $m=1,2\ldots$ 
are smaller than $K^2$.
\end{lemma}

\begin{proof}
The proof is similar to that in \cite{L}, Lemma 2.5 and uses Lemma 1.1
together with uniform local finiteness of $Y$.  More precisely we use the
fact that the number of $j$-simplexes in $Y$ at distance at most one from a
 $j$-simplex $\sigma$ can be bounded \emph{independently} of  $\sigma$, 
say
$\#\{\tau \in |Y|_j : d(\tau, \sigma) \leq 1\} \leq b$.
\emph {A fortiori} the same is true (with the same constant $a$) for $Y_m$ 
for all $m$. We now estimate the $\ell^2$ norm of $\Delta \kappa$ for a
cochain $\kappa = \sum_\sigma a_\sigma \sigma$ (having identified a
simplex $\sigma$ with the dual cochain). Now 
$$
\Delta \kappa =
\sum_\sigma \left ( \sum_\tau D(\sigma, \tau ) a_\tau \right ) \sigma
$$
so that
$$
\sum_\sigma \left ( \sum_\tau D(\sigma, \tau) a_\tau \right )^2 \leq 
\sum_\sigma \left ( \sum_{d(\sigma, \tau )\leq 1} D(\sigma, \tau )^2 \right
) \left ( \sum_{d(\sigma, \tau )\leq 1} a_\tau^2 \right )
\leq C^2 b \sum_\sigma \sum_{d(\sigma, \tau ) \leq 1} a_\tau^2,
$$
where we have used Lemma 1.1 and Cauchy-Schwartz inequality.  In the last
sum above, for every simplex $\sigma$, $a_{\sigma}^2$ appears at most
$b$ times.  This proves that $\|\Delta \kappa \|^2 
\leq C^2 b^2 \|\kappa \|^2$.
Identical estimate holds (with the same proof) for $\Delta^{(m)}$
which yields the lemma if we set $K=\sqrt{C b}$.
\end{proof}

Observe that $\Delta_j$ can be regarded as a matrix with entries in  
${\mathbb Z} [\pi]$, since
by definition, the coboundary operator $d_j$ is a matrix with entries in  
${\mathbb Z} [\pi]$, and  
so is its adjoint $d_j^*$ as it is equal to the simplicial boundary
operator.
There is a natural trace for matrices with entries in 
${\mathbb Z} [\pi]$, viz.
$$
\Tr_{{\mathbb Z} [\pi]}(A)= \sum_i \Tr_{{\mathcal U} [\pi]}(A_{i,i}).
$$

\begin{lemma} $\;$Let $\pi$ be an amenable group and 
let $p(\lambda) = \sum_{r=0}^d a_r \lambda^r$ be a polynomial.  Then,
$$
\Tr_{{\mathbb Z} [\pi]}(p(\Delta_j))=
\lim_{m\rightarrow\infty}\frac{1}{N_{m}}\;
\Tr_{{\mathbb C}}\Big(p\Big(\Delta_j^{(m)}\Big)\Big).
$$
\end{lemma}

\begin{proof}  First observe that if $\sigma\in |Y_m|_j$ is such that 
$d(\sigma , \partial Y_m) > k$, then Lemma 1.2 implies that
$$
D_j^k(\sigma, \sigma) = \left<\Delta_j^k \delta_\sigma, \delta_\sigma\right>
= \left<\Delta_j^{(m)k} \delta_\sigma, \delta_\sigma\right> = 
D_j^{(m)k}(\sigma, \sigma).
$$
By  (\ref{inv}) and (\ref{vNt}) 
$$
\Tr_{{\mathbb Z} [\pi]}(p(\Delta_j))= \frac{1}{N_m} \sum_{\sigma\in |Y_m|_j}
\left< p(\Delta_j)\sigma , \sigma \right>.
$$
Therefore we see that
$$
\left| \Tr_{{\mathbb Z} [\pi]}(p(\Delta_j)) - \frac{1}{N_{m}}\;
\Tr_{{\mathbb C}}\Big(p\Big(\Delta_j^{(m)}\Big)\Big)\right| \le
$$
$$
\frac{1}{N_{m}} \, \sum_{r=0}^d \,  |a_r|  
\sum_{
\begin{array}{lcl} & \sigma\in |Y_m|_j \\
 & d(\sigma, \partial Y_m) \leq d
\end{array}}
\, \left( D^r(\sigma, \sigma) + D^{(m)r}(\sigma, \sigma)\right).
$$
Using Lemma 1.2, we see that there is a positive constant $C$ such that
$$
\left| \Tr_{{\mathbb Z} [\pi]}(p(\Delta_j)) - \frac{1}{N_{m}}\;
\Tr_{{\mathbb C}}\Big(p\Big(\Delta_j^{(m)}\Big)\Big)\right| \le
2\, \frac{\dot{N}_{m,d}}{N_{m}} \, \sum_{r=0}^d \, |a_r| \, C^r.
$$
The proof of the lemma is completed by taking the limit as 
$m\rightarrow\infty$.
\end{proof}

We next recall the following abstract lemmata of L\"uck \cite{L}.

\begin{lemma} $\;$Let $p_{n}(\mu)$ be a sequence of polynomials 
such that for
the characteristic function of the interval $[0,\lambda]$, 
$\chi_{[0,\lambda]}
(\mu)$, and an appropriate real number $L$,
$$
\lim_{n\rightarrow\infty}p_{n}(\mu)=\chi_{[0,\lambda]}(\mu)\;\;\mbox{ and }
\;\;|p_{n}(\mu)|\leq L
$$
holds for each $\mu\in[0,||\Delta_j||^{2}]$.  Then
$$
\lim_{n\rightarrow\infty}\Tr_{{\mathbb Z} [\pi]}(p_{n}(\Delta_j))=F(\lambda).
$$\end{lemma}

\begin{lemma} $\;$Let $G:V\rightarrow W$ be a linear map of 
finite dimensional
Hilbert spaces $V$ and $W$.  Let $p(t)=\det(t-G^{*}G)$ be the characteristic
polynomial of $G^{*}G$.  Then $p(t)$ can be written as $p(t)=t^{k}q(t)$
where $q(t)$ is a polynomial with $q(0)\neq 0$.  Let $K$ be a real number,
$K\geq\max\{1,||G||\}$ and $C>0$ be a positive constant with 
$|q(0)|\geq C>0$.
Let $E(\lambda)$ be the number of eigenvalues
 $\mu$ of $G^{*}G$, counted with multiplicity, 
satisfying
$\mu\leq\lambda$.  Then for $0<\lambda<1$, 
the following
estimate is satisfied.
$$
\frac{ \{E(\lambda)\}- \{E(0)\}}{\dim_{{\mathbb C}}V}\leq
\frac{-\log C}{\dim_{{\mathbb C}}V(-\log\lambda)}+
\frac{\log K^{2}}{-\log\lambda}\;.
$$\end{lemma}

\begin{proof}[Proof of theorem 2.1] 
Fix $\lambda\geq 0$ and define for $n\geq 1$ a continuous function
$f_{n}:{\mathbb R}\rightarrow {\mathbb R}$ by
$$
f_{n}(\mu)=\left\{\begin{array}{lcl}
1+\frac{1}{n} & \mbox{ if } & \mu\leq\lambda\\[+7pt]
1+\frac{1}{n}-n(\mu-\lambda) & \mbox{ if } & 
\lambda\leq\mu\leq\lambda+\frac{1}{n} \\[+7pt]
\frac{1}{n}     & \mbox{ if } & \lambda+\frac{1}{n}\leq \mu
\end{array}\right.
$$
Then clearly $\chi_{[0,\lambda]}(\mu)<f_{n+1}(\mu)<f_{n}(\mu)$ and $f_{n}
(\mu)\rightarrow\chi_{[0,\lambda]}(\mu)$ as $n\rightarrow\infty$ for all
$\mu\in[0,\infty)$.  For each $n$, choose a polynomial $p_{n}$ such that
$\chi_{[0,\lambda]}(\mu)<p_{n}(\mu)<f_{n}(\mu)$ holds for all 
$\mu\in[0,K^{2}]$.
We can always find such a polynomial by a sufficiently close approximation of
$f_{n+1}$.  Hence
$$
\chi_{[0,\lambda]}(\mu)<p_{n}(\mu)<2
$$
and
$$
\lim_{n\rightarrow\infty}p_{n}(\mu)=\chi_{[0,\lambda]}(\mu)
$$
for all $\mu\in [0,K^{2}]$.  Recall that $E_{m}(\lambda)$ denotes the number
of eigenvalues $\mu$ of $\Delta_j^{(m)}$ satisfying $\mu\leq\lambda$
and counted with multiplicity.  Note that 
$||\Delta_j^{(m)} || \leq K^{2}$
by Lemma 2.2.
$$
\begin{array}{lcl}
\displaystyle\frac{1}{N_{m}}\;
\Tr_{{\mathbb C}}\big(p_{n}(\Delta_j^{(m)})\big)
&=&\displaystyle \frac{1}{N_{m}}\sum_{\mu\in [0,K^2]}p_{n}(\mu)\\[+12pt]
\displaystyle & =& 
\displaystyle\frac{ E_{m}(\lambda)}{N_{m}}+\frac{1}{N_{m}}\left\{
\sum_{\mu\in [0,\lambda ]}(p_{n}(\mu)-1)+\sum_{\mu\in (\lambda , 
\lambda + 1/n]}p_{n}(\mu)\right.\\[+12pt]   
\displaystyle& &\displaystyle \hspace*{.5in}\left.+\;\sum_{\mu\in (\lambda
+ 1/n, K^2]}p_{n}(\mu)\right\}
\end{array}
$$
Hence, we see that
\begin{equation} \label{A}
F_{m}(\lambda)=\frac{ E_{m}(\lambda)}{N_{m}}
\leq\frac{1}{N_{m}}\;\Tr_{{\mathbb C}}\big(p_{n}(\Delta_j^{(m)})\big).
\end{equation}
In addition,
$$
\begin{array}{lcl}
\displaystyle\frac{1}{N_{m}}\;
\Tr_{{\mathbb C}}\big(p_{n}(\Delta_j^{(m)})\big)& \leq &
\displaystyle\frac{ E_{m}(\lambda)}{N_{m}} 
+\;\frac{1}{N_{m}}\sup\{p_{n}(\mu)-1:\mu\in[0,\lambda]\}\; 
E_{m}(\lambda) \\[+16pt]
\displaystyle &+&\displaystyle\;\frac{1}{N_{m}}\sup\{p_{n}(\mu):
\mu\in[\lambda,\lambda+1/n]\}\;
 (E_{m}(\lambda+1/n)-E_{m}(\lambda)) \\[+16pt]
\displaystyle &+&\displaystyle\;\frac{1}{N_{m}}\sup\{p_{n}(\mu):
\mu\in[\lambda+1/n,\;K^{2}]\}\;
 (E_{m}(K^{2})-E_{m}(\lambda+1/n)) \\[+16pt]
\displaystyle &\leq &\displaystyle\frac{ E_{m}(\lambda)}{N_{m}}+
\frac{ E_{m}(\lambda)}{nN_{m}}+
\frac{(1+1/n) (E_{m}(\lambda+1/n)-E_{m}(\lambda))}{N_{m}} 
\\[+16pt]
\displaystyle & &\displaystyle\hspace*{.5in}+\;
\frac{(E_{m}(K^{2})-E_{m}(\lambda+1/n))}
{nN_{m}} \\[+16pt]
\displaystyle &\leq &\displaystyle
\frac{ E_{m}(\lambda+1/n)}{N_{m}}+\frac{1}{n}\;
\frac{ E_{m}(K^{2})}{N_{m}} \\[+16pt]
\displaystyle &\leq & \displaystyle F_{m}(\lambda+1/n)+\frac{a}{n}
\end{array}
$$
since $E_m(K^2)=\dim C^j(Y_m) \leq aN_m$ for a positive constant
$a$ independent
of $m$.
It follows that 
\begin{equation} \label{B}
\frac{1}{N_{m}}\;\Tr_{{\mathbb C}}\big(p_{n}(\Delta_j^{(m)})\big)\leq F_{m}
(\lambda+1/n)+\frac{a}{n}.
\end{equation}
Taking the limit inferior in (\ref{B}) and the limit superior in (\ref{A}), 
as $m\rightarrow\infty$, we get that
\begin{equation} \label{C}
{\overline{F}}(\lambda)\leq \Tr_{{\mathbb Z} [\pi]}\big(p_{n}(\Delta_j)\big)
\leq\mbox{\underline{$F$}}(\lambda+1/n)+\frac{a}{n}.
\end{equation}
Taking the limit as $n\rightarrow\infty$ in (\ref{C}) and 
using Theorem 2.4, we see that
$$
{\overline{F}}(\lambda)\leq F(\lambda)
\leq\mbox{\underline{$F$}}^{+}(\lambda).
$$
For all $\varepsilon>0$ we  have
$$
F(\lambda)\leq\mbox{\underline{$F$}}^{+}(\lambda)\leq\mbox{\underline{$F$}}
(\lambda+\varepsilon)\leq {\overline{F}}(\lambda+\varepsilon)
\leq F(\lambda+\varepsilon).
$$
Since $F$ is right continuous, we see that
$$
F(\lambda)={\overline{F}}^{+}(\lambda)=\mbox{\underline{$F$}}^{+}
(\lambda)
$$
proving the first part of theorem 2.1.

Next we apply theorem 2.5 to $\Delta_j^{(m)}$.  Let $p_{m}(t)$ denote the
characteristic polynomial of $\Delta_j^{(m)}$ and 
$p_{m}(t)=t^{r_{m}}q_{m}(t)$ where
$q_{m}(0)\neq 0$.   
The matrix describing
$\Delta_j^{(m)}$ has integer entries.  Hence $p_{m}$ is a polynomial 
with integer
coefficients and $|q_{m}(0)|\geq 1$.  By Lemma 2.2 and Theorem 2.5 there
are constants $K$ and $C=1$ independent of $m$, such that
$$
\frac{F_{m}(\lambda)-F_{m}(0)}{a}\leq\frac{\log K^{2}}
{-\log\lambda}
$$
That is,
\begin{equation}\label{D}
F_{m}(\lambda)\leq F_{m}(0)-\frac{a\log K^{2}}{\log\lambda}.
\end{equation}
Taking limit inferior in (\ref{D}) as $m\rightarrow\infty$ yields
$$
\mbox{\underline{$F$}}(\lambda)\leq\mbox{\underline{$F$}}(0)
-\frac{a\log K^{2}}{\log\lambda}.
$$
Passing to the limit as $\lambda\rightarrow +0$, we get that
$$
 \mbox{\underline{$F$}}(0)=\mbox{\underline{$F$}}^{+}(0) 
\qquad \qquad
\mbox{and } \qquad\qquad {\overline{F}}(0)={\overline{F}}^{+}(0).
$$
We have seen already that ${\overline{F}}^{+}(0)=F(0)=\mbox{\underline{$F$}}
(0)$, which proves part ii) of Theorem 2.1.  Since $\displaystyle-
\frac{a\log K^{2}}{\log\lambda}$ is right continuous in $\lambda$, 
$$
{\overline{F}}^{+}(\lambda)\leq F(0)-\frac{a\log K^{2}}{\log\lambda}.
$$
Hence part iii) of Theorem 2.1 is also proved.
\end{proof}

We will need the following lemma in the proof of Theorem 0.2 in 
the last section. We
follow the proof of Lemma 3.3.1 in \cite{L}.
\begin{lemma}
$$
\int_{0+}^{K^2}\left\{\frac{F(\lambda)-F(0)}{\lambda}\right\} d\lambda\le
\liminf_{m\to\infty} 
\int_{0+}^{K^2}\left\{\frac{F_m(\lambda)-F_m(0)}{\lambda}\right\} d\lambda
$$
\end{lemma}

\begin{proof}
By Theorem 2.1, and the monotone convergence theorem, one has

\begin{align*}
\int_{0+}^{K^2}\left\{\frac{F(\lambda)-F(0)}{\lambda}\right\}d\lambda
& = \int_{0+}^{K^2}\left\{\frac{\underline{F}(\lambda)-
\underline{F}(0)}{\lambda}\right\}d\lambda\\
& =  \int_{0+}^{K^2}\liminf_{m\to\infty}\left\{\frac{F_m(\lambda)-
F_m(0)}{\lambda}\right\}d\lambda \\
& =  \int_{0+}^{K^2}\lim_{m\to\infty}\left(\inf\left\{\frac{F_n(\lambda)-
F_n(0)}{\lambda}|n \ge m\right\}\right)d\lambda \\
& =  \lim_{m\to\infty}\int_{0+}^{K^2}\inf\left\{\frac{F_n(\lambda)-
F_n(0)}{\lambda}|n \ge m\right\}d\lambda \\
& \le \liminf_{m\to\infty} \int_{0+}^{K^2}\left\{\frac{F_m(\lambda)-
F_m(0)}{\lambda}\right\}d\lambda. 
\end{align*}

\end{proof}
\section{Proofs of the main theorems}

In this section, we will prove the Amenable Approximation Theorem (Theorem
0.1)
of the introduction. We will also prove some related spectral results.

\begin{proof}[Proof of Theorem 0.1 (Amenable Approximation Theorem)]
Observe that 

\begin{align*}
\frac{b^j(Y_{m})}{N_{m}} & = 
\frac{\dim_{\mathbb C}\Big(\ker(\Delta_j^{(m)})\Big)}{N_{m}} \\
&=F_{m}(0). 
\end{align*}

Also observe that

\begin{align*}
b_{(2)}^j(Y:\pi) & = \dim_{\pi}\Big(\ker(\Delta_j)\Big)\\
&=F(0). 
\end{align*}

Therefore Theorem $0.1$ follows from Theorem 2.1 after taking the 
limit as $m\to\infty$.
\end{proof}

 Suppose that $M$ is a compact Riemannian manifold and $\Omega^{j}_{(2)}
(\widetilde{M})$ denote the Hilbert space of square integrable $j$-forms on
a normal covering space $\widetilde{M}$, with transformation group $\pi$.  
The Laplacian ${\widetilde{\Delta}}_{j}:\Omega^{j}_{(2)}
(\widetilde{M})\rightarrow\Omega^{j}_{(2)}(\widetilde{M})$ is essentially
self-adjoint and has a spectral decomposition 
$\{{\widetilde P}_{j}(\lambda):\lambda\in
[0,\infty)\}$ where each ${\widetilde P}_{j}(\lambda)$ has finite 
von Neumann trace.
The associated von Neumann spectral density function, 
${\widetilde F}(\lambda)$ is
defined as
$$
{\widetilde F}:[0,\infty)\rightarrow [0,\infty),\;\;\;
{\widetilde F}(\lambda)=
\Tr_{{\mathcal  U}(\pi)}
({\widetilde P}_{j}(\lambda)).
$$
Note that ${\widetilde F}(0)=b_{(2)}^{j}(\widetilde{M}:\pi)$ and that 
the spectrum of
$\widetilde{\Delta}_{j}$ has a gap at zero if and only if there is a
$\lambda>0$ such that
$$
{\widetilde F}(\lambda)={\widetilde F}(0).
$$

Suppose that $\pi$ is an amenable group.  Fix a triangulation 
$X$ on $M$.  Then the normal cover $\widetilde{M}$ has an induced 
triangulation
$Y$.  Let $Y_{m}$ denote be a subcomplex of $Y$ such that 
$\big\{Y_{m}\big\}^{\infty}_{m=1}$ 
is a regular exhaustion of $Y$.   Let $\Delta^{(m)}_j:C^{j}(Y_{m},
{\mathbb C})\rightarrow C^{j}(Y_{m},{\mathbb C})$ denote the 
combinatorial Laplacian, and
let $E^{(m)}_j(\lambda)$ denote
the number of eigenvalues $\mu$ of $\Delta^{(m)}_j$ which are less than or
equal to $\lambda$.  Under the hypotheses above we prove the following.

\begin{theorem}[Gap criterion] $\;$The spectrum of ${\widetilde{\Delta}}_j$
has a gap at zero if and only if there is a $\lambda>0$ such that
$$
\lim_{m\rightarrow\infty}\;
\frac{E^{(m)}_j(\lambda)-E^{(m)}_j(0)}{N_{m}}=0.
$$\end{theorem}

\begin{proof}
Let $\Delta_j:C^{j}_{(2)}(Y)\rightarrow
C^{j}_{(2)}(Y)$ denote the combinatorial Laplacian acting on $L^2$ 
j-cochains on $Y$.  
Then by \cite{GS}, \cite{E}, the von Neumann spectral density function $F$
of the combinatorial Laplacian $\Delta_j$ and the von Neumann 
spectral density 
function $\widetilde F$ of the 
analytic Laplacian ${\widetilde{\Delta}}_j$ are dilatationally
equivalent, that is, there are constants $C>0$ and $\varepsilon>0$ 
independent
of $\lambda$ such that for all $\lambda\in(0,\varepsilon)$,
\begin{equation} \label{star}
F(C^{-1}\lambda)\leq {\widetilde F}(\lambda)\leq F(C\lambda).
\end{equation}
Observe that $\frac{E^{(m)}_j(\lambda)}{N_{m}} = F_m(\lambda)$.
Therefore the theorem also follows from Theorem 2.1.
\end{proof}

There is a standing conjecture that the Novikov-Shubin invariants of a 
closed manifold are positive (see \cite{E}, \cite{ES} and \cite{GS} 
for its definition). The next theorem gives evidence supporting
this conjecture, at least in the case of amenable fundamental groups.

\begin{theorem}[Spectral density estimate] $\;$There are constants $C>0$
and $\varepsilon>0$ independent of $\lambda$, such that for all $\lambda\in
(0,\varepsilon)$
$$
{\widetilde F}(\lambda)-{\widetilde F}(0)\leq\frac{C}{-\log(\lambda)}\;.
$$\end{theorem}

\begin{proof}
This follows from Theorem 2.1 and Theorem 3.1
since $\widetilde{\Delta}_j$ has a gap at zero if and only if 
$
{\widetilde F}_j(\lambda)={\widetilde F}_j(0)
$  for some $\lambda>0$. 
\end{proof}

\section{On the determinant class conjecture}

There is a standing conjecture that any normal
covering space of a finite simplicial complex is of determinant class.
Our interest in this conjecture stems from our work on $L^2$ 
torsion \cite{CFM}, \cite{BFKM}. The $L^2$ torsion is a well defined 
element in the determinant line of the reduced $L^2$ cohomology, whenever
the covering space is of determinant class. In this section, we use 
the results 
of section 2 to prove that any {\em amenable} normal
covering space of a finite simplicial complex is of determinant class.

Recall that a covering space $Y$ of a finite simplicial complex $X$
is said to be of {\em determinant class} if, for $0 \le j \le n,$
$$ - \infty < \int^1_{0^+} \log \lambda d F (\lambda),$$
where $F(\lambda)$ denotes the von Neumann 
spectral density function of the combinatorial Laplacian $\Delta_j$
as in Section 2. 

Suppose that $M$ is a compact Riemannian manifold and $\Omega^{j}_{(2)}
(\widetilde{M})$ denote the Hilbert space of square integrable $j$-forms on
a normal covering space $\widetilde{M}$, with transformation group $\pi$.  
The Laplacian ${\widetilde{\Delta}}_{j}:\Omega^{j}_{(2)}
(\widetilde{M})\rightarrow\Omega^{j}_{(2)}(\widetilde{M})$ is essentially
self-adjoint and the associated von Neumann spectral density function, 
${\widetilde F}(\lambda)$ is
defined as in section 3. 
Note that ${\widetilde F}(0)=b_{(2)}^{j}(\widetilde{M}:\pi)$
Then $\widetilde{M}$ is said to be of  {\em analytic-determinant class}, 
if, for $0 \le j \le n,$
$$ - \infty < \int^1_{0^+} \log \lambda d {\widetilde F} (\lambda),$$
where ${\widetilde F}(\lambda)$ denotes the von Neumann spectral density 
function
of the analytic Laplacian ${\widetilde{\Delta}}_{j}$ as above. By results 
of  
Gromov and Shubin \cite{GS},
the condition that $\widetilde{M}$ is of  analytic-determinant class
is independent of the choice of Riemannian metric on $M$.

Fix a triangulation 
$X$ on $M$.  Then the normal cover $\widetilde{M}$ has an induced 
triangulation
$Y$. Then $\widetilde{M}$ is said to be of  {\em combinatorial-determinant 
class}
if $Y$ is of determinant class. Using results of Efremov \cite{E}, and 
\cite{GS} one sees that 
the condition that $\widetilde{M}$ is of combinatorial-determinant class
is independent of the choice of triangulation on $M$.

Using again results of \cite{E} and \cite{GS}, one observes as in 
\cite{BFKM} 
that the combinatorial and analytic 
notions of determinant class coincide, that is 
$\widetilde{M}$ is of combinatorial-determinant class if and only if 
$\widetilde{M}$ is of analytic-determinant class.  The appendix of \cite{BFK}
contains a proof that every residually finite covering of a compact manifold 
is of determinant class. Their proof is based on L\"uck's 
approximation of von Neumann spectral density functions \cite{L}.  Since 
an analogous approximation holds in our setting (cf. Section 2), 
we can apply the argument of \cite{BFK} to prove Theorem 0.2.

\begin{proof}[Proof of Theorem 0.2 (Determinant Class Theorem)]
Recall that the \emph{normalized} spectral density functions 
$$
F_{m} (\lambda) = \frac 1 {N_m} E_j^{(m)} (\lambda)
$$
are right continuous.
Observe that $F_{m}(\lambda)$ are step functions and
denote by
${\det}' \Delta_j^{(m)}$ the modified determinant of $\Delta_j^{(m)}$, 
i.e. the product
of all {\em nonzero} eigenvalues of $\Delta_j^{(m)}$.  Let $a_{m, j}$
be the
smallest nonzero eigenvalue and $b_{m, j}$ the largest eigenvalue of
$\Delta_j^{(m)}$.  Then, for any $a$ and $b$, such that 
$0 < a < a_{m, j}$ and
$b > b_{m,j}$,
\begin{equation}\label{one}
\frac 1 {N_m} \log {\det}' \Delta_j^{(m)} = \int_a^b \log \lambda d F_{m}
(\lambda).
\end{equation}
Integration by parts transforms the Stieltjes integral 
$\int_a^b \log \lambda d F_{m}
(\lambda)$ as follows.
\begin{equation}\label{two}
\int_a^b \log \lambda d F_{m} (\lambda)  =  (\log b) \big( F_{m} (b)
- F_{m} (0) \big)  -  \int_a^b \frac {F_{m} (\lambda) - F_{m} (0)} \lambda d
\lambda. 
\end{equation}
As before, $F(\lambda)$ denotes the spectral density function of the operator
$\Delta_j$ for a fixed $j$.
Recall that $F(\lambda)$ is continuous to the right in $\lambda$.  Denote
by ${\det}'_\pi\Delta_j$ the modified Fuglede-Kadison determinant 
(cf. \cite{FK}) 
of $\Delta_j$, 
that is, the Fuglede-Kadison determinant of $\Delta_j$ restricted to the 
orthogonal
complement of its kernel. It is given by the
following Lebesgue-Stieltjes integral,
$$
\log {\det}^\prime_\pi \Delta_j = \int^{K^2}_{0^+} \log 
\lambda d F (\lambda) 
$$
with $K$ as in Lemma 2.2, i.e. $ ||  \Delta_j ||  < K^2$,
where $||\Delta_j||$ is the operator norm of $\Delta_j$.

Integrating by parts, one obtains
\begin{align} \label{three}
\log {\det}^\prime_\pi (\Delta_j) & =  \log K^2 \big( F(K^2) - F(0) \big) 
\nonumber\\
& +  \lim_{\epsilon \rightarrow 0^+} 
\Big\{(- \log \epsilon) \big( F (\epsilon)
-F(0) \big) - \int_\epsilon^{K^2} \frac {F (\lambda) - 
F (0)} \lambda d \lambda
\Big\}.
\end{align}
Using the fact that $  \liminf_{\epsilon \rightarrow 0^+} (- \log \epsilon)
\big( F (\epsilon) - F (0) \big) \ge 0$
(in fact, this limit exists and is zero)
and $\frac {F (\lambda) - F (0)}
\lambda \ge 0$ for $\lambda > 0,$ one sees that
\begin{equation}\label{four}
\log {\det}^\prime_\pi (\Delta_j) \ge 
( \log K^2) \big(F (K^2) - F(0) \big) - 
\int_{0^+}^{K^2}
\frac {F(\lambda) - F (0)} \lambda d \lambda.
\end{equation}

We now complete the proof of Theorem 0.2.
The main ingredient is the estimate of $\log {\det}'_\pi(\Delta_j)$ in terms
of
$\log {{\det}}^\prime \Delta_j^{(m)}$ combined with the fact that 
$\log {\det}'
\Delta_j^{(m)}\ge 0$ as the determinant $\det' \Delta_j^{(m)}$ is a 
positive integer.
By Lemma 2.2, there exists a positive number $K$, $1 \le K < \infty$, 
such that, for $m \ge 1$,
$$
||  \Delta_j^{(m)} ||  \le K^2 \quad {\text{and}}\quad ||  \Delta_j ||  
\le K^2.$$
By Lemma 2.6,
\begin{equation}\label{five}
\int_{0^+}^{K^2} \frac {F (\lambda) - F (0)} \lambda d \lambda \le
\liminf_{m \rightarrow \infty} \int_{0^+}^{K^2} \frac
{F_{m} (\lambda) - F_{m} (0)} \lambda d \lambda.
\end{equation}
Combining (\ref{one}) and (\ref{two}) with the inequalities 
$\log {\det}' \Delta_j^{(m)}
\ge 0$, we obtain
\begin{equation}\label{six}
\int_{0^+}^{K^2} \frac {F_{m} (\lambda) - F_{m} (0)} \lambda d \lambda
\le (\log K^2) \big(F_{m} (K^2) - F_{m} (0)\big).
\end{equation}
From (\ref{four}), (\ref{five}) and (\ref{six}), we conclude that
\begin{equation}\label{seven}
\log {\det}'_\pi \Delta_j \ge (\log K^2) \big( F(K^2) - F(0) \big) 
- \liminf_{m \rightarrow \infty}(\log K^2) \big( F_{m} (K^2) - F_{m}
(0) \big).
\end{equation}
Now Theorem 2.1 yields
$$
F (\lambda) = \lim_{\epsilon \rightarrow 0^+} \liminf_{m \rightarrow
\infty}
F_{m} (\lambda + \epsilon)
$$
and
$$
F (0) = \lim_{m \rightarrow \infty} F_{m} (0).
$$
The last two equalities combined with (\ref{seven}) imply that 
$\log {\det}'_\pi \Delta_j \ge 0$. 
Since this is true for all $j=0,1,\ldots,\dim Y$, 
$Y$ is of determinant class.
\end{proof}

\end{document}